\documentclass[reqno]{amsproc}

\usepackage{amstext,amsmath,amssymb,amsfonts}
\usepackage[latin1]{inputenc}
\usepackage{epsfig}
\usepackage{hyperref}
\usepackage{color}
\usepackage{tikz}
\usepackage[all,arc,dvips,arrow,curve,cmactex]{xy}

\usepackage{latexsym}

\newtheorem{definition}{Definition}
\newtheorem{theorem}{Theorem}

\newtheorem{proposition}{Proposition}
\newtheorem{example}{Example}

\newcommand{\bea}{\begin{eqnarray}}
\newcommand{\eea}{\end{eqnarray}}
\newcommand{\beq}{\begin{equation}}
\newcommand{\eeq}{\end{equation}}
\newcommand{\enn}{\nonumber \end{equation}}

 \newcommand{\cF}{\mathcal{F}}

\title[Delta-Matroids Handle Slides]{Finding classes of delta-matroids closed under  handle slides\footnote{ R. C. A. was supported by ISF Grant 1050/16.}}

\author{R\'emi Cocou Avohou}
\address[R.C.A.]{
Einstein Institute of Mathematics, The Hebrew University of Jerusalem, Giv'at Ram, Jerusalem, 91904, Israel, \&
ICMPA-UNESCO Chair, 072BP50, Cotonou, Rep. of Benin,  \& Ecole Normale Superieure, B.P 72, Natitingou, Benin}
\email{avohou.r.cocou@mail.huji.ac.il}

\begin{document}

\maketitle 

\begin{abstract}
In this work, we study the operations of handle slides introduced recently for delta-matroids by  Iain Moffatt and Eunice Mphako-Banda. We then prove that the class of binary delta-matroids is the only class of delta-matroids closed under handle slides.
\\

\noindent MSC(2010): 05C10, 57M15

\noindent Key words: matroids, delta-matroids, ribbon graphs, handle slides
\end{abstract}

\tableofcontents

\section{Introduction}

Introduced in 1935 as a generalization of graphs and linear independence in vector spaces  \cite{whitney}, matroids are important combinatorial structures. 
There are several approaches for studying matroids. Among those, we will be interested in the approach emphasizing bases, \cite{bouchet4} : here a matroid is defined as a set system $M = (V, \mathcal{B})$, where $V$ is called the ground set of $M$ and $\mathcal{B}\subset \mathcal{P}(V)$ the base set of $M$ satisfying the following exchange axiom: 
For $B, B'\in \mathcal{B}$ and $x\in B-B'$, there exists $y\in B'-B$ such that $B\Delta \{x, y\}\in \mathcal{B}$.

Delta-matroids constitute an interesting generalization of matroids and have been introduced relatively recently \cite{Boucher3}. The main idea in this setting is to replace the differences in the previous axiom by symmetric differences and the resulting axiom is called the symmetric exchange axiom. Many of the nice properties associated with matroids extend to delta-matroids. In particular, the connection between delta-matroids and embedded graphs generalizes the classical connection between matroids and abstract graphs.

There is another combinatorial notion that relates to delta-matroids and that we will discuss hereafter. Ribbon graphs arise naturally as neighborhoods of graphs embedded in surfaces \cite{bollo}. In \cite{phantom}, a particular operation on ribbon graphs called  ``handle slide'' which ``slides'' the end of one edge over an edge adjacent to it in the cyclic order at a vertex has been extended to the class of delta-matroid. The same paper proves that the class of binary delta-matroid is stable under handle slides. Then the authors of that reference asked a question: ``What classes of delta-matroids are closed under handle slides?''

The goal of this work is to answer this question. 

The paper is structured as follows. Section 2 sets up our definitions of delta-matroid and the handle slide operation. Section 3 investigates if there are other classes of delta matroid stable under handle slides. Theorem \ref{theo:closed} finally gives a negative answer to that question.

\section{Delta-matroids and the handle slide operations}

In this section, we give a quick review of some important results related to delta-matroids and handle slide operations \cite{bouchet1, bouchet2}. 

\begin{definition}[Delta-matroid]
A delta-matroid $D=(E,\mathcal{F})$ is a finite set $E$ called ground  set together with a nonempty collection of subsets of the ground set, $\mathcal{F}$ called feasible set, that  satisfies the symmetric exchange axiom i.e $\forall F_1, F_2\in \mathcal{F}$
$$x\in F_1\Delta F_2\Rightarrow \exists y\in F_1\Delta F_2,\  F_1\Delta \{x,y\}\in \mathcal{F}.$$
\end{definition}
A  set system is a pair $(E,\mathcal{F})$ of a finite set $E$ together with a nonempty collection $\mathcal{F}$ of subsets of $E$ and a delta-matroid is a set system satisfying the symmetric exchange axiom.

Consider a symmetric binary matrix $A=(a_{vw}: v, w\in E)$. Let $A[W]=(a_{vw}: v, w\in W)$ for $W\subseteq E$ 
%and $A[W_1, W_2]=(a_{vw}: v\in W_1, w\in W_2)$, for $W_1, W_2\subseteq E$. U
under the convention that $A[\emptyset]$ has an inverse, $D(A)=(E, \{W: A[W] \text{ has an inverse}\})$ is a delta-matroid.

\begin{definition}[Twist]\label{def:twist}
Let $D=(E, \cF)$ be a set system. For $A\subseteq E$, the twist of $D$ with respect to $A$, denoted by $D\star A$ is given by $(E, \{A\Delta X | X\in \cF\})$. The dual of $D$ written $D^*$ is equal to $D\star E$.
\end{definition}

\begin{definition}[Binary delta-matroid]\label{def:binarydeltamat}
A  delta-matroid $D=D(E, \cF)$ is said to be binary if there exists $F\in \cF$ and a symmetric binary matrix $A$ such that $D=D(A)\star F$.
\end{definition}

\begin{definition}[Handle slides \cite{phantom}]
\label{def:hs}
Let $D=(E, \cF)$ be a set system, and $a, b\in E$ with $a\neq b$. We define $D_{ab}$ to be the set system $(E, \cF_{ab})$ where
$$\cF_{ab}=\cF\Delta \big\{X\cup a | X\cup b\in\cF \text{ and } X\subseteq E-\{a, b\}\big\}.$$
We call the move taking $D$ to $D_{ab}$ a handle slide taking $a$ over $b$.
\end{definition}

\begin{example}

Here are a few specific examples for binary delta-matroids. If $A$ is the adjacency matrix of a simple graph, then singleton sets are
not feasible.
For example, let $$
  A =\begin{pmatrix}
0&1&0&0\\
1&0&1&1\\0&1&0&1\\0&1&1&0
\end{pmatrix}.
$$
The vertices of the graph in Figure \ref{fig:adj} should be labeled so that
vertex 1 is the vertex of degree~1 in the figure,
vertex 2 is the vertex of degree~3, vertices 3 and 4 are the two vertices
of degree 2.
Then vertex $i$ corresponds to  row/column $i$ of $A$.

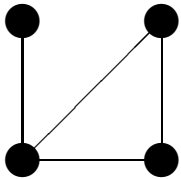
\begin{figure}[htb]\label{fig:adj}
\begin{center}
\unitlength=1.85mm
       \begin{picture}(16.25,16.25)
           \put(0.00,0.00){\line(1,0){10.00}}
           \put(0.00,0.00){\line(0,1){10.00}}
           \put(0.00,0.00){\line(1,1){10.00}}
           \put(10.00,0.00){\line(0,1){10.00}}
           \put(0.00,0.00){\circle*{2.50}}
           \put(10.00,0.00){\circle*{2.50}}
           \put(0.00,10.00){\circle*{2.50}}
           \put(10.00,10.00){\circle*{2.50}}

\end{picture}
\end{center}
   \caption{The graph of adjacency matrix $A$.}
\end{figure}
  The invertible
submatrices of $A$ correspond to the feasible sets $$\mathcal{F} = \{
\emptyset, \{1,2\}, \{2,3\}, \{2,4\}, \{3,4\}, \{1,2,3,4\} \}.$$
Here are some edge slides:
$$\mathcal{F}_{12} = \{ \emptyset, \{1,2\}, \{2,3\}, \{1,3\},\{1,4\}, \{2,4\},
\{1,2,3,4\} \},$$
$$\mathcal{F}_{21} = \{ \emptyset, \{1,2\}, \{2,3\}, \{2,4\}, \{3,4\}, \{1,2,3,4\} \},$$
$$\mathcal{F}_{23} = \{ \emptyset, \{1,2\}, \{2,3\}, \{3,4\}, \{1,2,3,4\} \},$$
$$\mathcal{F}_{32} = \{ \emptyset, \{1,2\}, \{2,3\}, \{1,3\}, \{2,4\}, \{1,2,3,4\} \}.$$

\end{example}

If a delta-matroid $D=(E, \cF)$ belongs to a class of delta-matroids closed under handle slides then for any $a, b\in E$, $D_{ab}$ is a delta-matroid and belongs to the same class. As consequence, if $D_{ab}$ is not a delta-matroid, this class can not exist. In general

\begin{definition}[Stable under handle slides]\label{def:classe}
Consider a (non-binary) delta-matroid $D=(E, \cF)$. We say that $D$ belongs to a class of (non-binary) delta-matroids closed or stable under handle slides if and only if for a sequence of handles slides, $(\cdots((D_{a_1b_1})_{a_2b_2})\cdots)_{a_nb_n}$ is a delta-matroid and belongs to the same class for any $a_1, b_1, \cdots, a_n, b_n\in E$ with $a_i\neq b_i$, $i=1, \cdots n$.
\end{definition}
A necessary condition for $D$ to belong to a class of delta-matroids closed under handle slides is that for any $a_1, b_1, \cdots, a_n, b_n\in E$, the set system $(\cdots((D_{a_1b_1})_{a_1b_2})\cdots)_{a_nb_n}$ is a delta-matroid.

\begin{definition}[Elementary minors] \label{def:minor}
Let $D=(E, \cF)$ be a delta-matroid. The elementary minors of $D$ at $e\in E$, are the delta-matroids $D-e$ and $D/e$ defined by:
$$D-e=\Big(E-e, \big\{F | F\subseteq E-e, F\in \cF\big\}\Big),$$
$$D/e=\Big(E-e, \big\{F | F\subseteq E-e, F\cup e\in \cF\big\}\Big).$$
The delta-matroid $D-e$ is called the deletion of $D$ along $e$, and $D/e$ the contraction of $D$ along $e$.
\end{definition}
A minor of a delta-matroid $D$ is obtained from $D$ by a sequence of deletions and contractions.

\begin{proposition}[Bouchet \cite{bouchet4}] \label{Prop:elemental}
If $D$ is a binary delta-matroid, then every elementary minor of $D$ is also a binary delta-matroid.
\end{proposition}

\begin{proposition}[Bouchet and Duchamp \cite{bouchet1}]\label{Prop:deltabinary}
A delta-matroid is binary if it has no minor isomorphic to a twist of $S_1$, $S_2$, $S_3$, $S_4$, or $S_5$, where
\bea
S_1&=&\Big(\{1, 2, 3\}, \big\{\emptyset, \{1, 2\}, \{1, 3\}, \{2, 3\}, \{1, 2, 3\}\big\}\Big),\cr\cr
S_2&=&\Big(\{1, 2, 3\}, \big\{\emptyset, \{1\}, \{2\}, \{3\}, \{1, 2\}, \{1, 3\}, \{2, 3\}\big\}\Big),\cr\cr
S_3&=&\Big(\{1, 2, 3\}, \big\{\emptyset, \{2\}, \{3\}, \{1, 2\}, \{1, 3\}, \{1, 2, 3\}\big\}\Big),\cr\cr
S_4&=&\Big(\{1, 2, 3, 4\}, \big\{\emptyset, \{1, 2\}, \{1, 3\}, \{1, 4\}, \{2, 3\}, \{2, 4\}, \{3, 4\}\big\}\Big),\cr\cr
S_5&=&\Big(\{1, 2, 3, 4\}, \big\{\emptyset, \{1, 2\}, \{1, 4\}, \{2, 3\}, \{3, 4\}, \{1, 2, 3, 4\}\big\}\Big).
\eea
\end{proposition}
\section{Finding classes of delta-matroids closed under handle slides}
This section investigates possible sub-classe(s) of non-binary delta-matroids closed under handle slides.

\begin{proposition}\label{Prop:minimalbinary}
None of the minimal non-binary delta-matroids $S_1$, $S_2$, $S_3$ or $S_5$ belongs to a class of delta-matroids stable under handle slides.
\end{proposition}
\proof
The proof of this lemma is direct by computing some handle slides of the delta-matroids $S_1$, $S_2$, $S_3$, $S_4$ and $S_5$. Let us compute $(S_1)_{12}$, $(S_2)_{12}$, $(S_3)_{23}$, $(S_4)_{12}$ and $(S_5)_{13}$. We obtain
\bea
(S_1)_{12}&=&\Big(\{1, 2, 3\}, \big\{\emptyset, \{1, 2\}, \{2, 3\}, \{1, 2, 3\}\big\}\Big),\cr
(S_2)_{12}&=&\Big(\{1, 2, 3\}, \big\{\emptyset, \{2\}, \{3\}, \{1, 2\}, \{2, 3\}\big\}\Big),\cr
(S_3)_{23}&=&\Big(\{1, 2, 3\}, \big\{\emptyset, \{3\}, \{1, 3\}, \{1, 2, 3\}\big\}\Big),\cr
(S_4)_{12}&=&\Big(\{1, 2, 3, 4\}, \big\{\emptyset, \{1, 2\}, \{2, 3\}, \{2, 4\}, \{3, 4\}\big\}\Big),\cr
(S_5)_{13}&=&\Big(\{1, 2, 3, 4\}, \big\{\emptyset, \{2, 3\}, \{3, 4\}, \{1, 2, 3, 4\}\big\}\Big),
\eea
which are not delta-matroids. The reason for this is that the symmetric exchange axiom is not satisfied for any of them. In conclusion, $S_1$, $S_2$, $S_3$, $S_4$ and $S_5$ do not belong to a class of delta-matroids stable under handle slides.
\qed

\begin{proposition}\label{Prop:stableunder}
Let $D=(E, \cF)$ be a delta-matroid and $a, b\in E$ ($a\neq b$). 
\begin{itemize}
\item[i] If $D_{ab}$ is a delta-matroid, then for any $e\in E-\{a, b\}$, $D_{ab}-e$ and $D_{ab}/e$ are also delta-matroids and
\bea
D_{ab}-e=(D-e)_{ab}, \quad D_{ab}/e=(D/e)_{ab}.
\eea
\item[ii] Let $a_1, b_1, \cdots, a_n, b_n\in E$. If the set system $(\cdots((D_{a_1b_1})_{a_2b_2})\cdots)_{a_nb_n}$ is a delta-matroid then 
$(\cdots(((D-e)_{a_1b_1})_{a_2b_2})\cdots)_{a_nb_n}$ and $(\cdots(((D/e)_{a_1b_1})_{a_2b_2})\cdots)_{a_nb_n}$ are also delta-matroids and,
in general, for every minor $T$ of $D$, $(\cdots((T_{a_1b_1})_{a_2b_2})\cdots)_{a_nb_n}$ is a delta-matroid.
\end{itemize}
\end{proposition}
\proof
Let us prove item i. The first assertion is a definition and the delta-matroids $D_{ab}-e$ and $D_{ab}/e$ are defined by: 
$$D_{ab}-e=(E-e, \cF_{a,b}-e), \quad D_{ab}/e=(E-e, \cF_{a,b}/e),$$
with $\cF_{a,b}-e=\{F | F\subseteq E-e, F\in \cF_{a,b}\}$ and $\cF_{a,b}/e= \{F-e | e\in F, F\in \cF_{ab}\}$. We need to prove that $\cF_{a,b}-e=(\cF-e)_{a,b}$ and $\cF_{a,b}/e=(\cF/e)_{a,b}$ where $(\cF-e)_{a,b}$ and $(\cF/e)_{a,b}$ are the feasible of $(D-e)_{ab}$ and $(D/e)_{ab}$ respectively. We have:
\bea
\cF_{a,b}-e&=&\Big\{F | F\subseteq E-e, F\in \cF \Delta \big\{X\cup a | X\cup b\in\cF \text{ and } X\subseteq E-\{a, b\}\big\}\Big\},\cr
&=&\big\{F | F\subseteq E-e, F\in \cF\big\}\Delta \big\{X\cup a | X\cup b\in\cF \text{ and } X\subseteq E-\{e, a, b\}\big\}\cr &=& (\cF-e)_{a,b},
\eea
and
\bea
\cF_{a,b}/e&=&\Big\{F | F\subseteq E-e, F\cup e\in \cF \Delta \big\{X\cup a | X\cup b\in\cF \text{ and } X\subseteq E-\{a, b\}\big\}\Big\},\cr
&=&\big\{F | F\subseteq E-e, F\cup e\in \cF\big\}\Delta \big\{X\cup a | X\cup b\in\cF \text{ and } X\subseteq E-\{e, a, b\}\big\}\cr &=& (\cF/e)_{a,b}. 
\eea
Item i follows.

We now concentrate on item ii. From Definition \ref{def:classe}, the delta-matroid $D$ belongs to a class of delta-matroids closed under handle slides if and only if for a sequence of handles slides, $(\cdots((D_{a_1b_1})_{a_1b_2})\cdots)_{a_nb_n}$ is a delta-matroid for any $a_1, b_1, \cdots, a_n, b_n\in E$ and belongs to the same class. From item i, we have $(\cdots(((D-e)_{a_1b_1})_{a_1b_2})\cdots)_{a_nb_n}=(\cdots((D_{a_1b_1})_{a_1b_2})\cdots)_{a_nb_n}-e$ and $(\cdots(((D/e)_{a_1b_1})_{a_1b_2})\cdots)_{a_nb_n}=(\cdots((D_{a_1b_1})_{a_1b_2})\cdots)_{a_nb_n}/e$. From this and the definition of a minor, $(\cdots((D_{a_1b_1})_{a_1b_2})\cdots)_{a_nb_n}$ is a delta-matroid. The result follows.
\qed
\begin{theorem}\label{theo:closed}
The only class of delta-matroids closed under handle slides is the class of binary delta-matroids.
\end{theorem}
\proof
We already know that the set of binary delta-matroids is closed under handle slides.  Let $D=(E, \cF)$ be a non-binary delta-matroid. If $D$ is minimal, Proposition \ref{Prop:minimalbinary} claims that $D$ does not belong to any class of delta-matroids closed under handle slides. Otherwise, from Proposition \ref{Prop:deltabinary}, $D$ has at least one minor $T$ isomorphic to a twist of $S_{i}$,  $i=1,\cdots, 5$. The fact that $(S_i)_{ab}$ is not a delta matroid does not imply necessarily that $(S_i\star A)_{ab}$ is not a delta-matroid. As example $(S_2\star \{1\})_{12}=\Big(\{1, 2, 3\}, \big\{\emptyset, \{1, 2\}, \{1, 3\}, \{2\}, \{3\}, \{1, 2, 3\}\big\}\Big)$ is a delta-matroid but $(S_2)_{12}$ is not. We then need to study all the possible cases:

$\bullet$ Assume that $T$ is isomorphic to a twist of $S_1$ i.e $T\cong S_1\star A$ with $A=$ $\emptyset$, $\{1\}$, $\{2\}$, $\{3\}$, $\{1, 2\}$, $\{1, 3\}$, $\{2, 3\}$, $\{1, 2, 3\}$. One has

\bea\label{eq:s1}
(S_1)_{12}&=&\Big(\{1, 2, 3\}, \big\{\emptyset, \{1, 2\}, \{2, 3\}, \{1, 2, 3\}\big\}\Big),\cr
(S_1\star \{1\})_{12}&=&\Big(\{1, 2, 3\}, \big\{\{2\}, \{3\}, \{1, 2, 3\}, \{1, 3\}, \{2, 3\}\big\}\Big),\cr
(S_1\star \{2\})_{12}&=&\Big(\{1, 2, 3\}, \big\{\{2\}, \{3\}, \{1, 2, 3\}, \{1, 3\}\big\}\Big),\cr
(S_1\star \{1, 3\})_{12}&=&\Big(\{1, 2, 3\}, \big\{\emptyset, \{1\}, \{2\}, \{1, 2\}, \{2, 3\}\big\}\Big),\cr
(S_1\star \{2, 3\})_{12}&=&\Big(\{1, 2, 3\}, \big\{\emptyset, \{1\}, \{1, 2\}, \{2, 3\}\big\}\Big),\cr
(S_1\star \{1, 2\})_{21}&=&(S_1)_{12}\star \{1, 2\},\cr
(S_1\star \{3\})_{12}&=&(S_1)_{12}\star \{3\},\cr
(S_1\star \{1, 2, 3\})_{21}&=&(S_1)_{12}\star \{1, 2, 3\}.
\eea
None of the sets system in \eqref{eq:s1} is a delta-matroid then from Proposition  \ref{Prop:stableunder}, $T\cong S_1\star A$, $A\subseteq \{1, 2, 3\}$ does not belong to any class of delta-matroids stable under handle slides.

$\bullet$ Let $T$ be isomorphic to a twist of $S_2$ i.e $T\cong S_2\star A$ with $A=$ $\emptyset$, $\{1\}$, $\{2\}$, $\{3\}$, $\{1, 2\}$, $\{1, 3\}$, $\{2, 3\}$, $\{1, 2, 3\}$. Remark that $(S_2)_{ab}$ ($a\neq b$) is not a delta-matroid for any $a, b=1, 2, 3$ ($a\neq b$) and then

\bea\label{eq:s2}
(S_2)_{12}&=&\Big(\{1, 2, 3\}, \big\{\emptyset, \{1\}, \{2\}, \{3\}, \{1, 2\}, \{2, 3\}\big\}\Big),\cr
(S_2\star \{1\})_{23}&=&(S_2)_{23}\star \{1\},\cr
(S_2\star \{2\})_{13}&=&(S_2)_{13}\star \{2\},\cr
(S_2\star \{3\})_{12}&=&(S_2)_{12}\star \{3\},\cr
(S_2\star \{1, 2\})_{21}&=&(S_2)_{12}\star \{1, 2\},\cr
((S_2\star \{1, 3\})_{23})_{12}&=&\Big(\{1, 2, 3\}, \big\{\emptyset, \{2\}, \{3\}, \{2, 3\},  \{1, 2, 3\}\big\}\Big)=((S_2\star \{2, 3\})_{13})_{12},\cr
(S_2\star \{1, 2, 3\})_{21}&=&(S_2)_{12}\star \{1, 2, 3\}.
\eea
The sets system in \eqref{eq:s2} are not delta-matroids then from Proposition \ref{Prop:stableunder}, $S_2\star A$, $A\subseteq \{1, 2, 3\}$ does not belong to any class of delta-matroids closed under handle slides and the same holds for $T$.  

$\bullet$ We now assume that $T$ is isomorphic to a twist of $S_3$ i.e $T\cong S_3\star A$ with $A=$ $\emptyset$, $\{1\}$, $\{2\}$, $\{3\}$, $\{1, 2\}$, $\{1, 3\}$, $\{2, 3\}$, $\{1, 2, 3\}$.

\bea\label{eq:s3}
(S_3)_{23}&=&\Big(\{1, 2, 3\}, \big\{\emptyset, \{3\}, \{1, 3\}, \{1, 2, 3\}\big\}\Big),\cr
(S_3\star \{1\})_{23}&=&(S_3)_{23}\star \{1\},\cr
(S_3\star \{2\})_{12}&=&\Big(\{1, 2, 3\}, \big\{\emptyset, \{2\}, \{2, 3\}, \{1, 2, 3\}\big\}\Big),\cr
(S_3\star \{3\})_{13}&=&\Big(\{1, 2, 3\}, \big\{\emptyset, \{3\}, \{2, 3\}, \{1, 2, 3\}\big\}\Big),\cr
(S_3\star \{1, 2\})_{13}&=&(S_3\star \{3\})_{13},\cr
(S_3\star \{1, 3\})_{12}&=&(S_3\star \{2\})_{12},\cr
(S_3\star \{2, 3\})_{32}&=&(S_3)_{23}\star \{2, 3\},\cr
(S_3\star \{1, 2, 3\})_{23}&=&(S_3)_{23}\star \{1, 2, 3\}.
\eea
In equation \eqref{eq:s3}, none of the sets system is a delta-matroids then from Proposition \ref{Prop:stableunder}, $S_3\star A$, $A\subseteq \{1, 2, 3\}$ does not belong to a class of delta-matroids closed under handle slides and the same holds for $T$.  

$\bullet$ Suppose that $T$ is isomorphic to a twist of $S_4$ i.e $T\cong S_4\star A$ with $A=$ $\emptyset$, $\{1\}$, $\{2\}$, $\{3\}$, $\{4\}$, $\{1, 2\}$, $\{1, 3\}$, $\{1, 4\}$, $\{2, 3\}$, $\{2, 4\}$, $\{3, 4\}$, $\{1, 2, 3\}$, $\{1, 2, 4\}$, $\{1, 3, 4\}$, $\{2, 3, 4\}$, $\{1, 2, 3, 4\}$. Remark that $(S_4)_{ab}$ ($a\neq b$) is not a delta-matroid for any $a, b=1, 2, 3, 4$ ($a\neq b$) and 

\bea
(S_4\star A)_{ab}=(S_4)_{ab}\star A, 
\eea
for $a, b\notin A$ and for $a, b\in A$, 
\bea
(S_4\star A)_{ba}=(S_4)_{ab}\star A.
\eea

Having these with  Proposition \ref{Prop:stableunder} implies that $S_4\star A$, $A\subseteq \{1, 2, 3, 4\}$ does not belong to a class of delta-matroids closed under handle slides and the same holds for $T$.  

Then from Proposition \ref{Prop:stableunder}, $T$ does not belong to any class of delta-matroids closed under handle slides and the result follows.

$\bullet$ We now study the lat case where $T$ is isomorphic to a twist of $S_5$ i.e $T\cong S_5\star A$ with $A=$ $\emptyset$, $\{1\}$, $\{2\}$, $\{3\}$, $\{4\}$, $\{1, 2\}$, $\{1, 3\}$, $\{1, 4\}$, $\{2, 3\}$, $\{2, 4\}$, $\{3, 4\}$, $\{1, 2, 3\}$, $\{1, 2, 4\}$, $\{1, 3, 4\}$, $\{2, 3, 4\}$, $\{1, 2, 3, 4\}$. One has

\bea\label{eq:s5}
(S_5)_{13}&=&\Big(\{1, 2, 3, 4\}, \big\{\emptyset, \{2, 3\}, \{3, 4\}, \{1, 2, 3, 4\}\big\}\Big),\cr
(S_5)_{24}&=&\Big(\{1, 2, 3, 4\}, \big\{\emptyset, \{1, 4\}, \{3, 4\}, \{1, 2, 3, 4\}\big\}),\cr
(S_5\star \{1, 2\})_{14}&=&\Big(\{1, 2, 3, 4\}, \big\{\emptyset, \{2, 4\}, \{3, 4\}, \{1, 2, 3, 4\}\big\}\Big)=(S_5\star \{3, 4\})_{14},\cr
(S_5\star \{1, 4\})_{12}&=&\Big(\{1, 2, 3\}, \big\{\emptyset, \{2, 4\}, \{2, 3\}, \{1, 2, 3, 4\}\big\}\Big)=(S_5\star \{2, 3\})_{12}.
\eea
The sets system in \eqref{eq:s5} are not delta-matroids and for $A\subset \{1, 2, 3, 4\}$ with $A\neq$ $\{1, 2\}$, $\{1, 4\}$, $\{2, 3\}$, $\{3, 4\}$, else $A$ contains none of the pairs $\{1, 3\}$ and $\{2, 4\}$ or contains one of them. Using the relations
\bea
(S_5\star A)_{ab}=(S_5)_{ab}\star A, 
\eea
for $a, b\notin A$ and 
\bea
(S_5\star A)_{ba}=(S_5)_{ab}\star A,
\eea
for $a, b\in A$, with  Proposition \ref{Prop:stableunder}, $S_5\star A$, $A\subseteq \{1, 2, 3, 4\}$ does not belong to a class of delta-matroids closed under handle slides and the same holds for $T$.  

In conclusion, the minor $T$ of $D$ does not belong to any class of delta-matroids stable under handle slides. This ends the proof of the theorem.
\qed

\section*{Acknowledgments}
I would like to thank Brigitte Servatius for her helpful comments and suggestions.

\vspace{0.5cm}

%\begin{center}
%\rule{3cm}{0.01cm}
%\end{center}

\end{document}